
\input amstex 
\documentstyle{amsppt} 

\magnification1200
\pagewidth{30pc} \pageheight{47pc} 


\define\N{\hbox{\bf N}}

\define\rd{\hbox{\rm d}} 

\define\Spec#1{\Phi _ {#1}}

\topmatter 
\title Strong Regularity for Uniform Algebras 
\endtitle 
\author J. F. Feinstein and D. W. B. Somerset\endauthor 
\leftheadtext{J. Feinstein and D. Somerset}%
\rightheadtext{Strong Regularity for Uniform Algebras}%

\address School of Mathematical Sciences, 
University of Nottingham, NG7 2RD, U.K. 
\endaddress 
\email Joel.Feinstein\@nottingham.ac.uk\endemail 

\address 
Department of Mathematical Sciences, 
University of Aberdeen, AB24 3UE, U.K. 
\endaddress 

\email ds\@maths.abdn.ac.uk\endemail 



\abstract 
A survey is given of the work on strong regularity for uniform 
algebras over the last thirty years, and some new results are 
proved, including the following. Let $A$ be a uniform algebra 
on a compact space $X$ and let $E$ be 
the set of all those points $x\in X$ such that $A$ is not strongly regular 
at $x$. If $E$ has no non-empty, perfect subsets then $A$ is normal, and $X$ 
is the character space of $A$. If $X$ is either $[0,1]$ or the circle $T$ and 
$E$ is meagre with no non-empty, perfect subsets then $A$ is trivial. These 
results extend Wilken's work from 1969. It is also shown 
that every separable Banach function algebra which has character space equal 
to either $[0,1]$ or $T$ 
and has a countably-generated ideal lattice is uniformly 
dense in the algebra of all 
continuous functions. 
\endabstract 

\thanks \endthanks 

\endtopmatter 

\document

The study of strong regularity for uniform 
algebras was initiated by Wilken in 1969 \cite{Wi2}. He had two main results. 
The first (a lemma) was that if a uniform algebra $A$ is strongly regular on a 
space $X$ then $X$ is necessarily the character space of $A$, and $A$ is 
normal. The second 
result was that there is no strongly regular uniform algebra on the unit 
interval $[0,1]$ other than the trivial one $C[0,1]$ itself. In this paper 
we survey the subsequent developments of the study, and extend 
Wilken's 
results. Various pertinent examples are given. 
\smallskip 
We begin by introducing the terminology, definitions and notation which we 
shall need, along with some background results concerning strong regularity. 
\smallskip 
\remark{Terminology and notation} In our terminology a 
{\it compact space\/} is a compact, Hausdorff topological space. 
For any compact space $X$ we denote the algebra of 
all continuous, complex-valued functions on $X$ by $C(X)$. 
\endremark 
\par 
Let $A$ be a commutative, 
unital Banach algebra. We denote by $\Spec{A}$ the character space of $A$. 
 
\definition{Definition} Let $X$ be a compact space. A {\it function 
algebra\/} 
on $X$ 
is a subalgebra of $C(X)$ which contains the constant functions and separates 
the points of $X$. A function algebra $A$ on $X$ is {\it trivial\/} if 
$A=C(X)$. 
A {\it Banach function algebra\/} on $X$ is a function algebra on $X$ with a 
complete algebra norm. 
A {\it uniform algebra\/} on $X$ is a Banach function algebra on $X$ whose 
norm is the uniform norm on $X$. 
\enddefinition 
\smallskip 
Every commutative, unital, semisimple Banach algebra may be regarded, via 
the Gelfand transform, as a Banach function algebra on its character space. 
In this paper we give results for unital Banach algebras. 
However, our results carry through without much difficulty to the non-unital 
case (with appropriate definitions) simply by considering the standard 
unitisation of any algebra under consideration. 
\remark{Notation} Let $A$ be a Banach function algebra on a compact 
space $X$, and let $x \in X$. We denote the evaluation character 
at $x$ by $\epsilon_x$. 
We define the ideals $J_x$, $M_x$ as follows: 
$$M_x = \{f \in A: f(x)=0\} = \ker \epsilon_x;$$ 
$$J_x = \{f \in A: f^{-1}(\{ 0 \}) \hbox{ is a neighbourhood of } x \}.$$ 
Let $K$ be a closed subset of $X$. Then we define the ideal $I(K)$ by 
$$I(K)=\{f \in A: f(K) \subseteq \{0\}\}.$$ 
\endremark 

It is standard to identify the set of evaluation characters 
$\{\epsilon_x: x\in X\}$ with $X$. 

\definition{Definition} Let $A$ be a Banach function algebra on a compact 
space 
$X$. Then $A$ is {\it regular on \/} $X$ if, for each closed set $F$ contained 
in 
$X$ and each $x\in X\backslash F$, there is an $f \in A$ with 
$f(F) \subseteq \{0\}$ and $f(x)=1$; 
$A$ is {\it regular \/} if it is regular on $\Spec{A}$. 
The algebra $A$ is {\it normal on \/} $X$ if, for every pair of disjoint 
closed 
sets 
$E$, $F$ 
contained in $X$, there is an $f \in A$ with $f(E) \subseteq \{ 0 \}$ and 
$f(F) \subseteq \{ 1 \}$; $A$ is {\it normal\/} if it is normal on 
$\Spec{A}$. Now let $x \in X$. Then $A$ is {\it strongly regular at\/} $x$ if 
$J_x$ is dense in $M_x$. The algebra $A$ is {\it strongly regular on $X$} if 
it is strongly regular at 
every point of $X$; $A$ is {\it strongly regular\/} if it is strongly regular 
on $\Spec{A}$. 
Finally, $A$ {\it has spectral synthesis\/} if $A$ is regular and if every 
closed ideal of $A$ is 
an intersection of maximal ideals. 
\enddefinition 
\smallskip 
A simple compactness argument shows that every strongly regular Banach 
function algebra is 
regular. It is standard, see e.g. \cite{St\rm, 27.2}, that every regular 
Banach 
function 
algebra is normal. 
\smallskip 
We now come to Wilken's results (mentioned earlier) about strong regularity 
for uniform algebras. The first is this. 
\smallskip 
\proclaim{Theorem 1 \cite{Wi2}} 
Let $A$ be a uniform algebra on a compact space $X$. Suppose that $A$ is 
strongly 
regular on $X$. Then 
\item{(a)} $\Spec{A}=X$; 
\item{(b)} $A$ is normal. 
\endproclaim 
\smallskip 
The second matter addressed by Wilken was the question of whether there 
actually are any non-trivial, strongly regular uniform algebras. He was 
able to show that there are no non-trivial, strongly regular uniform 
algebras on $[0,1]$. In 1973 his method was extended by Chalice \cite{Ch} 
to show that there are no non-trivial, strongly regular uniform algebras 
on the circle $T$ either. 
In 1974 some further results in this direction were 
obtained by Batikjan \cite{Ba} for locally connected, compact spaces 
$X$, under additional technical conditions. 

Also in \cite{Ch}, 
Chalice gave an example of a non-trivial uniform algebra which was 
strongly regular at all points of a dense subset of its character space. 
He raised the question whether a uniform algebra could be 
strongly regular at a non-peak point, and this question was resolved 
affirmatively by Wang in 1975 \cite{Wa}. The general question of the 
existence of non-trivial, strongly regular uniform algebras remained open, 
however. 

Meanwhile in 1987 Mortini gave an elementary proof of Wilken's result 
\cite{M}. 
Although not mentioned there, his method actually applies to all Banach 
function algebras. 
\smallskip 
\proclaim{Theorem 2} 
Let $A$ be a Banach function algebra on a compact space $X$, and let $K$ be a 
closed subset of $X$ such that $I(K)$ is the zero ideal. Suppose that $A$ is 
strongly regular at every point of $K$. Then 
\item{(a)} $\Spec{A}=X=K$; 
\item{(b)} $A$ is normal. 
\endproclaim 
\demo{Proof} Suppose, for contradiction, that there exists a character $\phi$ 
in 
$\Spec{A} \backslash K$. 
For every $x \in K$ we have that $\phi \neq \epsilon _ x$, and so $M_x$ is not 
a subset 
of $\ker \phi$. Since $J_x$ is dense in $M_x$, there must be a function 
$f_x$ in 
$J_x \backslash \ker \phi$. Set 
$$N_x = \{ y \in X: f_x(y)=0\}.$$ 
Then, by the compactness of $K$, 
there are finitely many elements $x_1, x_2, \dots x_n$ of $K$ such that 
$$K \subseteq \bigcup_{k=1}^n N_{x_k}.$$ 
Set $f = f_{x_1}f_{x_2}\cdots f_{x_n}$. Then $\phi(f) \neq 0$. 
But $f \in I(K)$, 
and so $f=0$. This contradiction proves (a), and (b) is now immediate. 
\quad\qed 
\enddemo 

In 1992 the first author answered Wilken's question affirmatively by 
exhibiting a non-trivial, strongly regular uniform algebra \cite{F1}, using 
Cole's systems of root extensions \cite{Co}. The same method also produces 
examples of 
non-trivial, normal uniform algebras which fail the condition of strong 
regularity at 
exactly one point. An example where the character space 
is metrizable can be found in a remark on page 300 of \cite{F1} (see also 
\cite{F1\rm, Theorem 5.3}). Taking finite direct sums, one can then obtain 
examples which fail strong regularity at any specified finite number of 
points. 
A similar method produces examples which fail strong regularity at countably 
many points. 

\example{Example 3} There exists a normal uniform algebra $B$ on a compact, 
metrizable space $Y$ 
such that the set of points of $Y$ at which $B$ is not strongly regular 
is countably infinite. 
\endexample 
\par 
\demo{Proof} 
Let $A$ be a normal uniform algebra on a compact metric space $X$ 
such that there is exactly one point of $X$ at which $A$ is not 
strongly regular \cite{F1}. 
Let $B$ be the standard 
unitisation of the $c_0$ direct sum of countably many copies of $A$ 
(so that $B$ is the algebra of all those sequences $(f_n) \subseteq A$ 
which converge uniformly on $X$ to some constant function). Let 
$Y$ be the one-point compactification of $X \times \N$, 
and regard $B$ as a uniform algebra on $Y$ in the obvious way. 
\par 
It is elementary to check that this example has the required 
properties. 
\quad\qed 
\enddemo 

It is natural to wonder if there is a connection between strong regularity 
and Gleason parts. A curious example is given in \cite{F2} of a strongly 
regular uniform algebra for which every Gleason part except one is a 
singleton, 
and the other part has exactly two points. 

All these examples are normal, and this brings us to the first main question 
of the paper. Can Wilken's first result be extended to the case where 
strong regularity holds except on some small exceptional set? Notice that 
Mortini's argument does not help with this question, as the following example 
shows. 

\example{Example 4} Let $X = \{ 0, 1, 1/2, 1/3, \dots\}$, and let $A$ be 
the restriction to 
$X$ of the disc algebra. By the identity principle, the restriction map is an 
isomorphism, 
and so $A$ is a Banach function algebra on $X$ (where the norm is the 
uniform norm of the functions on the closed unit disc). It is easy to see 
that this Banach function 
algebra is strongly regular at every point of $X$ except for the point $0$. 
But, of course, $\Spec{A} \neq X$, and $A$ is not normal. 
\endexample 

To extend Wilken's result therefore, we shall have return to his original 
method for uniform algebras. We shall need the following results. The first is 
due to Rudin from 1957. 
 
\par 
\proclaim{Theorem 5 \cite{R\rm, Theorem 4}} Let $X$ be a compact space which 
has no 
non-empty, perfect subsets. Then there are no non-trivial uniform algebras 
on $X$. 
\endproclaim 
\par 

The next is more recent. We need a couple of definitions. 

\definition{Definition} Let $X$ be a compact space, and let $A$ be a Banach 
function algebra on $X$. A point $x\in X$ is an {\it R-point \/} if $y\in X$ 
with $y\ne x$ implies that $M_y\not\supseteq J_x$. Evidently $x$ is an R-point 
whenever $A$ is strongly regular at $x$. A simple compactness argument 
shows that $A$ is regular on $X$ if and only if every point of $X$ is 
an R-point. 

The algebra $A$ is {\it $2$-local on X} if the following 
condition holds: if $f \in C(X)$ is such that 
there are elements $g_1$ and $g_2$ in $A$ so 
that every point of $X$ has neighbourhood on which $f$ agrees with 
either $g_1$ or $g_2$, then $f$ is in 
fact in $A$. 
Every local Banach function algebra is, of course, $2$-local on its 
character space. Thus every normal Banach function algebra is $2$-local on 
its character space. It is shown in \cite{Wi1\rm, 2.3, 3.1} that a 
uniform algebra whose character space is equal to either $[0,1]$ or $T$ is 
$2$-local 
on its character space. 

\enddefinition 

\proclaim{Theorem 6 \cite{FS2}} 
Let $A$ be a Banach function algebra on a compact space $X$, and let $y\in X$. 
The set $F=\{ x\in X: M_y\supseteq J_x\}$ is closed, and is connected if $A$ 
is 
$2$-local on $X$. If $A$ is a uniform algebra on $X$ then $F$ is either a 
singleton or 
contains a non-empty perfect subset. 
\endproclaim 

Since points of $F\backslash\{ y\}$ are not R-points, it follows that 
the set of non-R-points has a non-empty, perfect subset if $A$ is a 
uniform algebra which is not regular on $X$. 
\smallskip 
We also recall some standard facts about Jensen measures. Wilken used the 
existence of these in his proof of Theorem 1. 

\definition{Definition} Let $A$ be a uniform algebra on a 
compact space $X$, and let 
$\phi \in \Spec{A}$. Then a {\it Jensen} measure for $\phi$ is a regular, 
Borel probability measure $\mu$ on $X$ such that, for all $f \in A$, 
$$\log |\phi(f)| \leq \int_X{\log |f(x)| \rd \mu(x)}\eqno(1.1)$$ 
(where $\log(0)$ is defined to be $-\infty$). 
\enddefinition 
\smallskip 
It is standard (see, for example, \cite{G\rm, p.33}) that every $\phi \in 
\Spec{A}$ 
has a Jensen measure supported on $X$, and that each such measure 
{\it represents} $\phi$, i.e., for all $f \in A$, 
$$\phi(f) = \int_X{f(x) \rd \mu(x)}.$$ 

We are now ready for our first new result. 

\proclaim{Theorem 7} Let $A$ be a uniform algebra on a compact space $X$. 
Set $$E=\{x\in X:\overline{J_x} \neq M_x\}.$$ 
Suppose that $E$ has no non-empty, perfect subsets. Then 
\item{(a)} $\Spec{A} = X$; 
\item{(b)} $A$ is normal. 
\endproclaim 
\par 
\demo{Proof} For (a), let $\phi \in \Spec{A}$. Let $\mu$ be a Jensen measure 
for $\phi$ 
supported on $X$, and let $F$ be the closed support of $\mu$. 
Suppose first that $F \backslash E \neq \emptyset$. Let $x\in F\backslash E$. 
Then it follows from (1.1) that $J_x \subseteq M_{\phi}$. 
Since $\overline{J_x}=M_x$, we obtain 
$M_x \subseteq M_\phi$, and so $\phi = \epsilon_x$. 
Now suppose instead that $F \subseteq E$. Set $K = F \cup \{\phi\}$. Then $K$ 
is a compact subset of $\Spec{A}$. Since $F \subseteq E$ it follows that 
$F$ has no non-empty, perfect subsets. The same must also be true for $K$. 
Thus, by Theorem 5, $A|K$ is uniformly dense in $C(K)$. On the 
other hand, $\delta_\phi - \mu$ is a regular Borel measure on $K$ 
which annihilates $A$. Thus we must have $\delta_\phi=\mu$. Since $\mu$ is 
supported on $X$, it follows that $\phi \in X$, as required. 

For (b), observe that if $A$ were not normal on $X$ then by Theorem 6 there 
would be a non-empty, perfect subset of $X$ consisting of points which were 
not R-points, and a fortiori not points of strong regularity. Thus $E$ would 
have a non-empty, perfect subset. This contradiction establishes that $A$ is 
normal. 
\quad\qed 
\enddemo 
\par 
 
The following examples illustrate what can happen once a uniform algebra fails 
to be strongly regular on a large enough subset. 

\example{Examples 8} 
(a) Let $S=\bigcup_{n=1}^{\infty}\left\{\left( 
1-\frac1{n}\right) 
e^{{2\pi ki\over n}}: 1\le k\le n\right\}$, and let $X=S\cup T$. Let $A$ be 
the restriction to 
$X$ of the disc algebra. By the maximum modulus principle, the 
restriction map is an isometric isomorphism, 
so $A$ is a uniform algebra on $X$. Every point of $S$ is isolated, so $A$ is 
strongly regular at every point of $S$, which is a dense open subset of $X$. 
But, of course, $\Spec{A} \neq X$, and $A$ is not normal. 

(b) Let $A$ be the uniform algebra of continuous functions on 
a solid cylinder which are analytic on the base of the cylinder (often 
called `the tomato-can algebra'). Strong regularity holds everywhere 
except the base, hence on a 
dense open subset of $\Spec{A}$, but $A$ is not normal. 

(c) Let $A$ be the uniform algebra obtained by restricting $H^{\infty}$ 
(the algebra 
of bounded functions on the disc which are analytic on the open disc) to 
the fibre of its maximal ideal space associated with a point on the unit 
circle, see \cite{H\rm, p.187ff}. Then $A$ is regular on its Shilov boundary, 
but $A$ is not normal. 
This example shows that the requirement 
of strong regularity outside the exceptional set in Theorem 7 
cannot be relaxed to regularity (i.e. simply requiring that 
each $x$ outside the exceptional set should be an R-point). In this 
example the exceptional set of non-R-points is actually empty. 

(d) O'Farrell \cite{O} has given an example of a normal 
uniform algebra $R(X)$, consisting of the closure of the 
rational functions having poles off the Swiss cheese $X$, with the 
property that there are continuous point derivations on a set 
of positive measure. Since strong regularity fails for normal uniform 
algebras wherever there 
is a continuous point derivation, it follows that $R(X)$ fails 
strong regularity at uncountably many points. 

\endexample 

We turn now to consider the second part of Wilken's work on strong regularity 
for uniform algebras, namely that every strongly regular uniform algebra on 
$[0,1]$ is trivial. The question of whether there are any non-trivial 
uniform algebras with character space equal to $[0,1]$ seems still to be 
open. With Theorem 6, however, we are able to push Wilken's method a little 
further. 

\definition{Definition} A closed set $E \subseteq X$ is a {\it peak set} for 
$A$ if there is a function $f$ in $A$ such that $f$ is constantly $1$ 
on $E$, but such that $|f(t)|<1$ for all $t \in X\backslash E$. 
A {\it peak point} for $A$ is a point $x \in X$ such that the 
set $\{x\}$ is a peak set for $A$. 

Now let $A$ be a uniform algebra on a compact space $X$. 
It is standard that when $X$ is metrizable (which occurs if and 
only if $A$ is separable) the set of peak points is a dense, $G_{\delta}$ 
in the Shilov boundary of $A$, see II.11.2 and II.12.10 of 
\cite{G}. Furthermore countable intersections 
and finite unions of peak sets are again peak sets for $A$, see Section 
II.12 of \cite{G}. In particular the union of two peak points is a peak set. 
\enddefinition 

For the next lemma, let us say that an ideal $I$ in a Banach algebra 
{\it factors} if, for all $f$ in $I$, there are $g$, $h$ in $I$ 
such that $f=gh$. By Cohen's factorization theorem, see 
\cite{P\rm, 5.2.2}, $I$ factors 
whenever $I$ has a bounded approximate identity. If $x$ is a peak point 
for a uniform algebra 
then $M_x$ has a bounded approximate identity, and hence factors. 

\proclaim{Lemma 9} Let $A$ be a Banach function algebra on a compact space 
$X$. Let 
$F$ be the set of all those $x$ in $X$ such that 
$M_x$ factors and $A$ is strongly regular at $x$. Suppose 
that $x\in F$, and that $y\in X$ with $A$ strongly regular at $y$. 
Then $J_x \cap J_y$ is dense in 
$M_x \cap M_y$. 
\endproclaim 
\par 
\demo{Proof} Let $f \in M_x \cap M_y$. Because $M_x$ factors we can 
write $f=gh$ where $g,h\in M_x$. Since $f\in M_y$, at least one ($h$ 
say) of $g$ and $h$ belongs to $M_y$. Choose sequences 
$(g_n) \subseteq J_x$ and $(h_n) \subseteq J_y$ converging to 
$g$ and $h$ respectively. Then $(g_n h_n) \subseteq J_x \cap J_y$ and 
this sequence converges to $f$. The result follows.\quad\qed 
\enddemo 

\proclaim{Proposition 10} Let $X$ be $[0,1]$ or $T$, and let $A$ be a uniform 
algebra on 
$X$. Suppose that $A$ is $2$-local on $X$ and that $X$ has a dense 
subset consisting of peak points at which $A$ is 
strongly regular. Then $A$ is trivial. 
\endproclaim 
\par 
\demo{Proof} The proofs for $[0,1]$ and for $T$ are subtly different, due 
to the existence of end-points in $[0,1]$. Let $S$ be the dense 
subset of $X$ consisting of peak points of $X$ at which 
$A$ is strongly regular. 

First suppose that $X=[0,1]$. 
Let $x$ be in $S$. We proceed as in Wilken's original 
proof: choose $f \in A$ such that $f(x)=1$ and 
such that $|f(t)|<1$ for all other points of $[0,1]$. 
Choose a sequence $g_n$ in $J_x$ converging to $1-f$ in $A$. 
Define $h_n(t)$ to be $1$ on $[0,x]$ and $1-g_n(t)$ for all 
other $t$. Then the functions $h_n$ agree locally on $X$ 
with either $1$ or $1-g_n$, so each $h_n$ is in $A$, since $A$ is 
$2$-local on $X$. 

Also, the functions $h_n$ converge to a function $h$ in $A$ 
which is $1$ on $[0,x]$, but with $|h(t)|<1$ for all 
other $t$. Thus $[0,x]$ is a peak set for $A$, and, similarly, 
so too is $[x,1]$. Since $S$ is dense, and countable intersections and finite 
unions 
of peak sets are peak sets, we see that $[x,y]$ 
is a peak set for all 
$x$, $y$ in $[0,1]$, and thus that every closed subset of 
$[0,1]$ is a peak set. 
It follows, by \cite{Br\rm, 2.4.3}, that for every closed subset 
$E$ of $[0,1]$ the restriction of $A$ to $E$ is closed in $C(E)$. 
Thus $A=C([0,1])$ by Glicksberg's theorem, see \cite{St\rm, 13.5}. 

Now suppose instead that $X=T$. Let $E$ be any closed arc in $T$ 
with endpoints in $S$. We shall show that $E$ is a peak set for 
$A$. The result will then follow, as before. 
Let $x$, $y$ be the end-points of $E$. Then (by 
the earlier remarks) $\{x,y\}$ is a peak set for $A$. 
Choose $f \in A$ such that $f(x)=f(y)=1$ and such that $|f(t)|<1$ 
for all other points of $T$. By Lemma 9, we can find a 
sequence $(g_n)$ in $J_x \cap J_y$ converging to $1-f$ in $A$. 

Define $h_n(t)$ to be $1$ for $t \in E$, and set $h_n(t)=1-g_n(t)$ 
for all other $t$. As above, each $h_n$ is in $A$. The 
functions $h_n$ converge uniformly to a function $h$ which is 
constantly $1$ on $E$, but such that $|h(t)|<1$ for all other 
$t$. Thus $E$ is a peak set for $A$, as claimed. 
The rest of the proof is identical to the proof above for the 
interval, noting, of course, that $T$ itself is trivially a peak set 
for $A$. 
\quad\qed 
\enddemo 

\proclaim{Theorem 11} Let $X$ be $[0,1]$ or $T$, and let $A$ be a uniform 
algebra on 
$X$. Set $$E=\{x\in X:\overline{J_x} \neq M_x\}.$$ 
Suppose that $A$ is $2$-local on $X$ and that $E$ is meagre. Then 
$A$ is trivial. 
\endproclaim 

\demo{Proof} Since $A$ is $2$-local on $X$, Theorem 6 shows that if $A$ were 
not regular on $X$ then there would be an interval consisting of non-R-points, 
hence points of $E$. This would contradict the hypothesis that $E$ is meagre. 
Hence $A$ is regular on $X$, so $X$ is the Shilov boundary of $A$. This 
implies 
that the set of peak points 
is a dense $G_{\delta}$ in $X$, so its intersection with the complement of the 
meagre set $E$ gives a dense subset of peak points at which 
strong regularity holds. Thus $A$ is trivial by Proposition 10. \quad\qed 
\enddemo 

If $A$ is a uniform algebra with $\Spec{A}$ equal to either $[0,1]$ 
or $T$ then $A$ is $2$-local on $X=\Spec{A}$ \cite{Wi1\rm, 2.3, 3.1}. Thus 
if the set $E$ above 
is meagre then $A$ is trivial. 

\proclaim{Corollary 12} Let $A$ be a uniform algebra on $X$, where $X$ is 
either $[0,1]$ or $T$. Set $$E=\{x\in X:\overline{J_x} \neq M_x\}.$$ 
Suppose that $E$ is meagre and has no non-empty, perfect subsets. Then 
$A$ is trivial. 
\endproclaim 

\demo{Proof} Theorem 7 shows that $X$ is the character space of $A$, and 
that $A$ is normal on $X$, and hence $2$-local. The result now follows from 
Theorem 11. \quad\qed 
\enddemo 

In particular, if the exceptional set $E$ is countable then $A$ is trivial. 
The following example shows that things can go wrong if the exceptional 
set is permitted to have a non-empty perfect subset. 

\example{Example 13} Let $A$ be the non-trivial uniform algebra on the 
Cantor set described in Theorem 9.3 of \cite{We}. The character space 
of $A$ is the whole of the Riemann sphere (see 9.2 and 9.3$'$ of \cite{We}), 
so \cite{St\rm, 27.3} shows that $A$ is not normal on the Cantor set. 
Let $B$ be the algebra of continuous functions on $[0,1]$ whose 
restrictions to the Cantor set lie in $A$. Then $B$ is a non-normal 
uniform algebra on $[0,1]$, strongly regular on a dense, open subset of 
$[0,1]$. The character space of $B$ is not equal to $[0,1]$ since it contains 
a copy of the 
Riemann sphere. Theorem 11 shows that $B$ is not $2$-local on $[0,1]$. 
\endexample 

For the last part of the paper we consider 
a condition on Banach algebras which has been of interest recently, and which 
has consequences for strong regularity. 
For a Banach algebra $A$, let $Id(A)$ be the lattice of closed, 
two-sided ideals of $A$. Then $Id(A)$ is {\it countably-generated \/} if there 
is a countable subset $S$ of $Id(A)$ such that for all $I\in Id(A)$ 
$$I=\left(\bigcup \{ J\in S: J\subseteq I\}\right)^-.$$ 
Examples of Banach algebras with countably-generated ideal lattices include 
separable C$^*$-algebras, TAF-algebras, and separable 
Banach algebras with spectral synthesis, see \cite{So}. For further results 
on this property, see \cite {Be1} and \cite{Be2}. We shall need the following 
result of Beckhoff's. 

\proclaim{Lemma 14} \cite {Be1} Let $A$ be a separable Banach algebra. The 
following are 
equivalent. 
\item{(a)} $Id(A)$ is countably-generated. 
\item{(b)} There is a countable subset $B$ of $A$ such that 
for all $I\in Id(A)$, $I=\overline{I\cap B}$. 
\endproclaim 

\smallskip 

\proclaim{Proposition 15} Let $A$ be a separable Banach function algebra on a 
compact space $X$. Set $$E=\{x\in X:\overline{J_x} \neq M_x\}.$$ If $Id(A)$ is 
countably-generated then $E$ is a meagre subset of $X$. 
\endproclaim 
\par 
\demo{Proof} By Lemma 14 there is a countable subset $B$ of $A$ such that 
$I=\overline{I\cap B}$ for each $I\in Id(A)$. Note that for $b\in B$ and $x\in 
X$, a necessary condition for $b\in M_x\setminus \overline{J_x}$ is that $x$ 
should belong to the boundary of the zero set of $b$ (for otherwise either 
$b\notin M_x$ or $b\in J_x$). The boundary of the zero set of $b$ is a closed 
set without interior, hence meagre. Thus the set $$F=\{x\in X: M_x\cap B \neq 
\overline{J_x}\cap B\}$$ is meagre, being a countable union of meagre sets. 
But for $x\in X$, $x\notin F$ if and only if $M_x\cap B=\overline{J_x}\cap B$, 
which holds if and only if $M_x=\overline{J_x}$. 
Thus $E=F$, so $E$ is meagre. 
\quad\qed 
\enddemo 

\example{Examples 16} (a) For $\alpha\in (0,1)$, let $A=lip_{\alpha}[0,1]$ be 
the little 
Lipschitz algebra on the metric space $([0,1], d^{\alpha})$, where 
$d^{\alpha}(x,y)=|x-y|^{\alpha}$, $(x,y\in [0,1])$. 
Then $A$ is separable and $\Spec{A}=[0,1]$. It was shown in \cite{Sh\rm, 
Corollary 4.3} that $A$ has spectral synthesis. Hence $Id(A)$ is 
countably-generated, by the remarks above. 

(b) Let $A$ be the Banach function algebra consisting of 
those continuous functions $f$ on $[0,1]$ which are differentiable 
at $0$, with norm given by $$\Vert f\Vert=\Vert f\Vert_{\infty}+\sup\left\{ 
\left\vert{f(t)-f(0)\over t}\right\vert:t\in (0,1]\right\},$$ where $\Vert 
\cdot\Vert_{\infty}$ denotes the uniform norm. Then $A$ is separable and 
normal on $[0,1]$ and strong regularity fails only at the point $0$. The ideal 
lattice of $A$ is countably-generated, see \cite{Be2\rm, p.455}. 

(c) Let $A=C^1[0,1]$. Then $A$ is separable and normal but there is no 
point in $[0,1]$ at which $A$ is strongly regular. Hence $Id(A)$ is not 
countably-generated. 

(d) Let $A$ be the disc algebra. Then $A$ is separable but there is no 
point of the disc at which $A$ is strongly regular. Hence again $Id(A)$ is not 
countably-generated. 
\endexample 

The next lemma shows that strong regularity is preserved 
under uniform closure. 

\proclaim{Lemma 17} Let $B$ be a Banach function algebra on a compact space 
$X$ and let $x \in X$. Let $A$ be the uniform closure of $B$ in $C(X)$. 
If $B$ is strongly regular at $x$, then $J_x\cap B$ is uniformly dense in 
$M_x$ 
(where $M_x$ and $J_x$ denote the ideals of $A$). 
\endproclaim 

\demo{Proof} 
Let $f \in A$ with $f(x)=0$. Choose a sequence of functions 
$(f_n) \subseteq B$ converging uniformly to $f$ on $X$. By subtracting 
the constant $f_n(x)$ from $f_n$ if necessary, we may assume that $f_n(x)=0$. 
Then, since $B$ is strongly regular at $x$, there are functions $g_n$ in $B$ 
each of which vanishes on a neighbourhood of $x$, and such that 
$\|g_n-f_n\|_B < 1/n$. But the norm in $B$ must dominate the uniform norm on 
$X$ 
(which is at most equal to the spectral radius), so it follows that the 
sequence 
$(g_n)$ converges to $f$ uniformly on $X$. The result follows.\quad\qed 
\enddemo 

\proclaim{Theorem 18} Let $B$ be a separable Banach function algebra on $X$, 
where $X$ is either 
$[0,1]$ or $T$. Suppose either that $B$ is $2$-local on $X$, 
or that $\Spec{B}=X$. If $Id(B)$ is countably-generated, then $B$ is uniformly 
dense in $C(X)$. 
\endproclaim 

\demo{Proof} Let $A$ be the uniform closure of $B$ in $C(X)$. Suppose first 
that $B$ is 
$2$-local on $X$. As in Theorem 11, it follows from Theorem 6 and Proposition 
15 that $B$ must 
be regular on $X$. Hence $A$ is also regular on $X$, so $X$ is the Shilov 
boundary of $A$. Thus 
the set of peak points of $A$ is a dense 
$G_{\delta}$ of $X$. Proposition 15 shows that strong regularity holds for $B$ 
on a dense $G_{\delta}$ of $X$, and strong 
regularity holds for $A$ at every 
point of $X$ where it held for $B$, by Lemma 17. Thus $X$ has a dense subset 
$S$ consisting of 
peak points of $A$ where strong regularity holds for $A$. We now follow the 
proof of Proposition 10. 

First suppose that $X=[0,1]$. 
Let $x$ be in $S$. Choose $f \in A$ such that $f(x)=1$ and 
such that $|f(t)|<1$ for all other points of $[0,1]$. 
Lemma 17 shows that there is a sequence $(g_n)$ of elements of $B$ 
each vanishing in a neighbourhood of $x$, and converging in $A$ to $1-f$. 
Define $h_n(t)$ to be $1$ on $[0,x]$ and $1-g_n(t)$ for all 
other $t$. Then each $h_n$ agrees locally on $X$ 
with either $1$ or $1-g_n$, so each $h_n$ is in $B$, since $B$ is $2$-local. 
The rest of the proof for $X=[0,1]$ now follows as in Proposition 10. 

Now suppose that $X=T$. Let $E$ be any closed arc in $T$ 
with endpoints $x,y\in S$. Choose $f \in A$ such that $f(x)=f(y)=1$ 
and such that $|f(t)|<1$ 
for all other points of $T$. Then $(1-f)\in M_x\cap M_y$ (where again we use 
$M_x$, 
$J_x$, etc. to denote ideals of $A$), so as in Lemma 9 there exist $p\in 
M_x$ and $q\in M_y$ such that $(1-f)=pq$. By Lemma 17 there are sequences 
$(p_n)$ in $J_x\cap B$ and $(q_n)$ in $J_y\cap B$ converging to $p$ and $q$ 
respectively in $A$. Hence if $g_n=p_nq_n$ then $(g_n)\subseteq J_x\cap 
J_y\cap B$ and 
$(g_n)$ converges to $1-f$ in $A$. 
Define $h_n(t)$ to be $1$ for $t \in E$, and $h_n(t)=1-g_n(t)$ 
for all other $t$. As above, each $h_n$ is in $B$. The proof is now 
concluded as in Proposition 10. 

Finally, suppose that $\Spec{B}=X$. Then $\Spec{A}=X$ too, and again 
Proposition 15 and Lemma 17 shows that strong regularity holds for $A$ 
on a dense $G_{\delta}$ of $X$. Hence $A$ is trivial, by the remark after 
Theorem 11. 
\quad\qed 
\enddemo 

Of course, the conditions of Theorem 18 are by no means necessary 
for the conclusion, as Example 16(c) shows, for instance. 
\smallskip 
Let us conclude by mentioning a couple of open questions. It follows from 
Theorem 5.1 of \cite{Wh} 
that whenever $A$ is a normal 
uniform algebra on a compact space $X$ and $x\in X$ is such that $J_x$ 
has a bounded approximate identity, then $A$ is strongly regular at $x$. 
It appears to be unknown, however, whether 
a normal uniform algebra must be strongly 
regular at each of its peak points (recall that if a point $x$ is a peak point 
then 
$M_x$ has a bounded approximate identity). 
If this should be the case, then 
the methods above will show immediately that every normal uniform algebra 
on $[0,1]$ or $T$ is trivial. For some positive results in this 
direction, see \cite{Ba}. 

It is also unknown whether or not there are any non-trivial 
uniform algebras that have spectral synthesis. 
In \cite{FS1} there is an example 
of a strongly regular uniform algebra $A$ such that every maximal ideal of 
$A$ has a bounded approximate identity, but such that $A$ does not have 
spectral synthesis. In the same paper it is shown that the 
method of taking systems of root extensions 
cannot produce a uniform algebra with spectral synthesis unless it 
is applied to such an algebra in the first place.

\refstyle{A}

\enddocument